\documentclass[notitlepage,leqno,11pt]{article}
\usepackage{amssymb}
\catcode`\@=11 \@addtoreset{equation}{section}

\catcode`\@=12

\usepackage{latexsym}
\usepackage{textcomp}
\usepackage{amsmath}
\usepackage{amsfonts}
\usepackage{amssymb}
\usepackage{mathrsfs}

\renewcommand{\b }{\beta }
\renewcommand{\d}{\delta }

\newcommand{\pa}{{\partial}}
\newcommand{\D }{\Delta }

\newcommand{\e }{\varepsilon }

\renewcommand{\l }{\lambda }
\renewcommand{\L }{\Lambda }

\newcommand{\n }{\nabla }

\newcommand{\s }{\sigma }
\newcommand{\Sig }{\Sigma}

\renewcommand{\th }{\theta }
\renewcommand{\o }{\omega }
\renewcommand{\O }{\Omega }

\newcommand{\be}{\begin{equation}}
\newcommand{\ee}{\end{equation}}

\newcommand{\R}{\mathbb{R}}
\newcommand{\N}{\mathbb{N}}
\renewcommand{\S}{\mathbb{S}}

\newcommand{\C}{\mathbb{C}}

\newcommand{\de}{\partial}

\newcommand{\ti}{\tilde}

\newcommand{\ra}{{\rangle}}
\newcommand{\la}{{\langle}}

\newcommand{\calC }{\mathcal{C}}
\newcommand{\CS }{\mathcal{C}_{\Sigma}}

\newcommand{\calD }{\mathcal{D}}

\newtheorem{Theorem}{Theorem}[section]
\newtheorem{Lemma}[Theorem]{Lemma}
\newtheorem{Proposition}[Theorem]{Proposition}
\newtheorem{Corollary}[Theorem]{Corollary}
\newtheorem{Remark}[Theorem]{Remark}

\def\proof{\noindent{{\bf Proof. }}}
\def\square{\vbox{
    \hrule height .4pt
    \hbox{\vrule width .4pt height 7pt \kern 7pt
       \vrule width .4pt}
    \hrule height .4pt }}

\def\QED{\hfill {$\square$}\goodbreak \medskip}

\def\R{{\mathbb R}}
\def\Rnk{{{\mathbb R}^{N}}}
\def\C{{\mathcal C}}
\def\S{{\mathbb S}}

\def\div{{\rm div}}

\font\sc=cmcsc9  \textwidth=14truecm
\author{Mouhamed Moustapha Fall
\footnote{\footnotesize{Universit\'e Catholique de Louvain-La-Neuve,
d\'epartement de  math\'ematique.
 Chemin du Cyclotron 2,  1348 Louvain-la-Neuve, Belgique.
E-mail: \textit{ mouhamed.fall@uclouvain.be},
\textit{mouhamed.m.fall@gmail.com}. } }}

\begin{document}

\title{A note on Hardy's inequalities with  boundary
singularities}

\date{}
\maketitle

\bigskip

\noindent {\footnotesize{\bf Abstract.}  Let $\O$ be a smooth
bounded domain in $\R^N$ with   $N\ge 1$. In this paper we study the
Hardy-Poincar\'e inequalities with  weight function singular at the
boundary of $\O$. In particular we give sufficient  conditions so
that the best constant is achieved.}
\bigskip\bigskip

\section{Introduction}

Let $\Omega$ be a  domain in $\R^N$, $N\ge 1$, with $0\in\de\O$ and
$p>1$ a real number.
 In this note, we are interested in finding minima to the following
 quotient
\begin{equation}
\label{eq:problem-m} \mu_{\l,p}(\O):= \inf_{u\in W^{1,p}_{0}(\O)}
~\frac{\displaystyle\int_{\O}|\nabla u|^p~dx-\l\int_{\O}|u|^p~dx}
{\displaystyle\int_{\O}|x|^{-p}|u|^p~dx}~,
\end{equation}
in terms of   $\lambda\in\R$ and $\O$.  If $\lambda=0$, we have the
{\em $\Omega$-Hardy constant}
\begin{equation}
\label{eq:mu} \mu_{0,p}(\O)= \inf_{u\in W^{1,p}_{0}(\O)}
~\frac{\displaystyle\int_{\O}|\nabla u|^p~dx}
{\displaystyle\int_{\O}|x|^{-p}|u|^p~dx}~
\end{equation}
which is the best constant in the Hardy inequality for maps
supported by $\Omega$.
The existence of extremals for  $\mu_{\l,2}(\O)$ was studied in
\cite{FaMu} while for $\mu_{0,2}(\O)$, one can see for instance
\cite{CaMuUMI}, \cite{CaMuPRSE}, \cite{PT} and \cite{NaC} for $\mu_{0,N}(\O)$.\\
Given a unit vector $\nu$ of $\R^N$, we consider the half-space
$H:=\{x\in\R^N\,:\, x\cdot\nu\geq0\}$. For $N=1$, the following
Hardy inequality is well known
\begin{equation}\label{eq:Hardy-1d}
 \left(\frac{p-1}{p}\right)^p\int_0^\infty t^{-p}|u|^p\,dt\leq
 \int_0^\infty|u'|^p\,dt\quad\forall u\in W^{1,p}_0(0,\infty).
\end{equation}
Moreover $\mu_{0,p}(H)=\left(\frac{p-1}{p}\right)^p$ is the
$H$-Hardy constant and it is not achieved, see \cite{OK} for
historical comments also.\\
For $N\geq2$, it was recently proved by  Nazarov \cite{Na} that the
$H$-Hardy constant is not achieved and
\be\label{eq:mpeq}
\mu_{0,p}(H):=\inf_{V\in
W^{1,p}_{0}(\S_+^{N-1})}\frac{\displaystyle\int_{\S_+^{N-1}}\left(\left(\frac{N-p}{p}
\right)^2|V|^2+|\n_\s V|^2 \right)^{\frac{p}{2}}d\s}{
\displaystyle\int_{\S_+^{N-1}}|V|^pd\s},
\ee
 where $\S_+^{N-1}$ is an
$(N-1)$-dimensional hemisphere. Notice that this problem always has
a minimizer by the compact embedding $L^p(\S_+^{N-1})
\hookrightarrow W^{1,p}_{0}(\S_+^{N-1})$. The quantity
$\mu_{0,p}(H)$ is explicitly known only in some special cases.
Indeed, $\mu_{0,2}(H)=\frac{N^2}{4}$ while for $p=N$ then
$\mu_{N,N}(H)$ is the first Dirichlet eigenvalue of the operator
$-\div(|\n u|^{N-2}\n u)$ in $W^{1,N}_0(\S_+^{N-1})$ with the
standard metric.\\
 Problem
(\ref{eq:problem-m}) carries some similarities with the questions
studied by Brezis and Marcus in \cite{BM}, where the weight is the
inverse-square of  the distance from the boundary of $\Omega$ and
$p=2$. We also deal with this problem in the present paper for all
$p>1$ in Appendix \ref{app:Hd}. We generalize here the existence
result obtained by R.Musina and the author in \cite{FaMu} for any
$p>1$ and $N\geq1$.
%
%
\begin{Theorem}\label{th:attained}
Let $p>1$ and $\O$ be a smooth bounded domain in $\R^N$, $N\geq1,$
with $0\in\de\O$. There exits $\l^*(p,\O) \in[-\infty,+\infty)$ such
that
\be\label{eq:mlssmp}
\mu_{\l,p}(\O)<\mu_{0,p}(H),\quad \forall \l>\l^*(p,\O).
\ee
The infinimum in \eqref{eq:problem-m} is attained for any
$\l>\l^*(p,\O)$.
\end{Theorem}
The existence of $\l^*(p,\O)$ comes from the fact that
 $$\sup_{\l\in\R}\mu_{\l,p}(\O)= \mu_{0,p}(H),$$
 see Lemma \ref{lem:supmp}. Now observe that the mapping $\l\mapsto\mu_{\l,p}$ is
non-increasing. Moreover, for bounded domains $\O$, letting $\l_1$
be the first Dirichlet eigenvalue of the $p$-Laplace operator
$-\div(|\n u|^{p-2}\n u)$ in $W^{1,p}_0(\O)$, it is plain that
$\mu_{\l_1,p}(\O)=0$. Then  we define
 $$
 \l^*(p,\O):=\inf\{\l\in\R\,:\, \mu_{\l,p}(\O)<\mu_{0,p}(H) \}
$$
so that $ \mu_{\l,p}<\mu_{0,p}(H)$ for all $\l>\l^*(p,\O)$. In
particular  $\l^*(p,\O)\leq\l_1$. On the other hand there are
various bounded smooth domains $\O$ with $0\in\de\O$ such that
$\l^*(p,\O)\in[-\infty,0)$, see Proposition \ref{prop:uppls} and
Proposition \ref{prop:uppls1}. Furthermore if $N=1$ then
$\mu_{0,p}(\R\setminus\{0\})=\left(\frac{p-1}{p}\right)^p=\mu_{0,p}(H)$
thus $\l^*(p,\O)\geq0$.

It is obvious that if $\O$ is contained in a half-ball centered at
the origin then $\mu_{0,p}(\O)=\mu_{0,p}(H)$ thus $\l^*(p,\O)\geq0$
and in addition
$$
\l^*(p,\O)=\inf_{u\in W^{1,p}_0(\O)}\frac{\displaystyle\int_{\O}|\n
u|^p\,dx- \mu_{0,p}(H)\displaystyle\int_{\O}|x|^{-p}|
u|^p\,dx}{\displaystyle\int_{\O}| u|^p\,dx}.
$$
We have obtained the following result.
\begin{Theorem}\label{th:lowbls} If $\O$ is contained in a
half-ball centered at the origin then there exists a  constant
$c(N,p)>0$ such that
\be\label{eq:lblst}
\l^*(p,\O)\geq\frac{c(N,p)}{\textrm{diam}(\O)^p}.
\ee
\end{Theorem}
%
The constant  $c(N,p)$ appearing  in \eqref{eq:lblst}  has the
property that $c(N,2)$ is the first Dirichlet eigenvalue of $-\D$ in
the unit disc of $\R^2$.
This type of estimates was first proved by Brezis-V\`{a}zquez in
\cite{BV} when $p=2$, $N\geq2$ and later on, extended to the case
$1<p<N$ by Gazzola-Grunau-Mitidieri in \cite{GGM} when dealing with
$\mu_{0,p}(\R^N\setminus\{0\}):=\left|\frac{N-p}{p}\right|^p$. More
precisely they proved the existence of a positive constant $C(N,p)$
such  that for any open subset $\O$ of $\R^N$, there holds
\be\label{eq:HardRN} \int_\Omega|\nabla
u|^p-\mu_{0,p}\left(\R^N\setminus\{0\}\right)\int_\Omega|x|^{-p}|u|^p\ge
C(N,p)\left(
\frac{\o_N}{|\Omega|}\right)^\frac{p}{N}\int_\Omega|u|^p~\!\quad\forall
u\in W^{1,p}_0(\O), \ee
where $|\O|$ is  the measure of $\O$ and $\o_N$ the measure of the
unit ball of $\R^N$. The constant $C(N,p)$ was explicitly given  and
$C(N,2)=c(N,2)$ as was obtained in \cite{BV}. The main ingredients
to prove \eqref{eq:HardRN} is the Schwarz symmetrization and a
"dimension reduction" via the transformation $x\mapsto
\frac{u}{\o}$, where $\o(x)=|x|^{\frac{p-N}{p}}$ satisfies
$$
\div(|\n
\o|^{p-2}\n\o)+\mu_{0,p}\left(\R^N\setminus\{0\}\right)\,|x|^{-p}\o^{p-1}=0\quad\textrm{in
}\R^{N}\setminus\{0\}.
$$
For $p=2$, the lower bound in \eqref{eq:lblst} was obtained in
\cite{FaMu} by a similar transformation and using the Poincar\'e
inequality on  $\S^{N-1}_+$. However, in view of \eqref{eq:mpeq},
such argument do not apply here when $p\neq2$ and $p\neq N$. By
analogy, to reduce the dimension, we will consider the mapping
$x\mapsto \frac{u}{v}$, where
$v(x):=|x|^{\frac{p-N}{p}}V\left(\frac{x}{|x|}\right)$ is a weak
solution to the  equation
$$
\div(|\n v|^{p-2}\n v)+\mu_{0,p}(H)\,|x|^{-p}|v|^{p-2}v=0\quad
\textrm{ in } \calD'(H)
$$
whenever $V$ is  a minimizer of \eqref{eq:mpeq}.  Then exploiting
the strict  convexity of the mapping  $a\mapsto |a|^p$,   estimate
\eqref{eq:lblst}, for $p\geq2$, follows immediately while the case
$p\in(1,2)$
    carries further difficulties as it can be seen in  Section \ref{ss:rt}. \\

 The argument to prove  the attainability  of $\mu_{\l,p}(\O)$ is taken
from de Valeriola-Willem \cite{dVW}. It
 allows to show that, up to a subsequence,
 the gradient of the Palais-Smale sequences  converges point-wise almost every where.
Therefore an application of the Brezis-Lieb lemma with some simples
arguments yields the existence of extremals.
\section{Hardy inequality with one point singularity}\label{s:H1ps}
Let $\calC$ be a proper cone in $\R^N$, $N\geq2$ and put
$\Sig:=\calC\cap\S^{N-1}$. It was shown in \cite{Na} that the
$\calC$-Hardy constant  is not achieved and it is given by
\begin{equation}\label{eq:mpH}
\mu_{0,p}(\calC)=\inf_{V\in
W^{1,p}_{0}(\Sig)}\frac{\displaystyle\int_\Sig\left(\left(\frac{N-p}{p}
\right)^2|V|^2+|\n_\s V|^2 \right)^{\frac{p}{2}}d\s}{
\displaystyle\int_{\Sig}|V|^pd\s}.
\end{equation}
 Letting $V\in W^{1,p}_{0}(\Sig)$ be the positive minimizer to this
quotient then the function
\be\label{eq:Hardeigf}
 v(x):=|x|^{\frac{p-N}{p}} V\left(\frac{x}{|x|}\right)
\ee
 satisfies
\be\label{eq:vweaksol}
\int_{\calC}|\n v|^{p-2}\n v\cdot \n
h=\mu_{0,p}(\calC)\,\int_{\calC}|x|^{-p} v^{p-1}h\quad\forall h\in
C^{1}_c(\calC).
\ee
%
Notice that
$\mu_{0,2}(\calC)=\left(\frac{N-2}{2}\right)^2+\l_1(\Sig)$, where
$\l_1(\Sig)$ is the first Dirichlet eigenfunction of the Laplace
operator on $\Sig$ endowed with the standard metric on $\S^{N-1}$.
This was obtained in \cite{PT}, \cite{NaC} and \cite{FaMu}.

%
\subsection{Existence}\label{ss:exist}
In this Section we show that the condition $\mu_{\l,p}(\Omega)<
\mu_{0,p}(H)$ is sufficient to guaranty the existence of a minimizer
for
$\mu_{\lambda,p}(\Omega)$.\\
 We emphasize  that throughout this section, $\O$ can to be taken to be  an
 open set satisfying the uniform sphere condition at $0\in\de\O$. Namely
there are balls $B_+\subset\O$ and $B_-\subset\R^N\setminus\O$ such
that  $\de{B_+}\cap\de{B_-}=\{0\}$. This holds if $\de\O$ is of
class $C^2$ at 0, see [\cite{GT} 14.6 Appendix]. We start with the
following approximate local Hardy inequality.
\begin{Lemma}\label{lem:loc-Hard}
Let $\O$ be a smooth  domain  in $\R^N$, $N\geq1,$ with $0\in\de\O$
and let $p>1$. Then for any $\e>0$ there exits $r_\e>0$ such that
\begin{equation}\label{eq:loc-H}
\mu_{0,p}(\O\cap B_{r_\e}(0))\geq \mu_{0,p}(H)-\e,
\end{equation}
where $B_{r}(0)$ is a ball of radius $r$  centered at 0.
\end{Lemma}
\proof If $N=1$ then \eqref{eq:loc-H} is an immediate consequence of
\eqref{eq:Hardy-1d}. From now on we can assume that $N\geq2$. We
denote by $N_{\de\O}$ the unit normal vector-field on $\de\O$. Up to
a rotation, we can assume that $N_{\de\O}(0)=E_N$, so that the
tangent plane of $\de\O$ at 0 coincides with
$\R^{N-1}=\textrm{span}\{E_1,\dots,E_{N-1}\}$. Denote by
$B_r^+=\{y\in B_r(0)\,:\, y^N>0\}$. For $r>0$ small, we introduce
the following system of coordinates centered at 0 (see
\cite{MMF-PJM}) via  the mapping $F:B_r^+\to \O$ given by
$$
 F(y)=\textrm{Exp}_0(\ti{y})+y^N\,
N_{\de\O}\left(\textrm{Exp}_0(\ti{y}) \right),
$$
where $\ti{y}=(y^1,\dots,y^{N-1})$ and  $ \ti{y}\mapsto
\textrm{Exp}_0(\ti{y})\in\de\O$ is the exponential mapping of
$\de\O$ endowed with the metric induced by $\R^N$. This coordinates
induces a metric on $\R^N$ given by
$g_{ij}(y)=\la\pa_iF(y),\pa_jF(y)\ra$ for $i,j=1,\dots,N$. Let $u\in
C^\infty_c(F(B_r^+))$ and put $v(y)=u(F(y))$ then
\begin{equation}\label{eq:coordvu}
\int_{F(B_r^+)}|\n u|^p\,dx=\int_{B_r^+}|\n
v|^p_g\sqrt{|g|}\,dy,\quad\int_{F(B_r^+)}|x|^{-p}|u|^p\,dx=\int_{B_r^+}|F(y)|^{-p}|v|^p\sqrt{|g|}\,dy,
\end{equation}
with $|g|$ stands for the determinant of the  $g$ while $|\n
v|^p_g=g(\n v,\n v)^{\frac{p}{2}}$. Since $|F(y)|=|y|+O(|y|^2)$ and
$g_{ij}(y)=\d_{ij}+O(|y|)$, we infer that
$$
\displaystyle\frac{\displaystyle  \int_{B_r^+}|\n
v|^p_g\sqrt{|g|}\,dy  }{\displaystyle
\int_{B_r^+}|F(y)|^{-p}|v|^p\sqrt{|g|}\,dy }\geq
(1-Cr)\displaystyle\frac{ \displaystyle\int_{B_r^+}|\n v|^p\,dy
}{\displaystyle \int_{B_r^+}|y|^{-p}|v|^p\,dy },
$$
for some  constant $C>0$ depending only on $\O$ and $p$. Furthermore
since $\mu_{0,p}(B_r^+)\geq\mu_{0,p}(H)$, using \eqref{eq:coordvu}
we conclude that
$$\mu_{0,p}(F(B_r^+))\geq(1-Cr)\mu_{0,p}(H).$$
 \QED
 We are in position to prove \eqref{eq:mlssmp} in  the following
\begin{Lemma}\label{lem:supmp}
Let $\O$ be a smooth  domain  in $\R^N$, $N\geq1,$ with $0\in\de\O$
and let $p>1$. Then there exists $\l^*(p,\O)\in [-\infty,+\infty)$
such that
$$
\mu_{\l,p}(\O)<\mu_{0,p}(H)\quad \forall \l>\l^*(p,\O).
$$
\end{Lemma}
\proof We first show that
\be\label{eq:supmlmp}
\sup_{\l\in\R}\mu_{\l,p}(\O)=\mu_{0,p}(H).
\ee
\textbf{Step 1}: We claim that
$\sup_{\l\in\R}\mu_{\l,p}(\O)\geq\mu_{0,p}(H)$.\\
%
For $r>0$ small, we  let $\psi\in C^\infty(B_{r}(0))$ with
$0\leq\psi\leq 1$, $\psi\equiv 0$ in $\R^N\setminus
B_{\frac{r}{2}}(0)$ and $\psi\equiv 1$ in $B_{\frac{r}{4}}(0)$.
For a fixed $\e>0$ small, there holds
\begin{eqnarray*}
\int_\O|x|^{-p} |u|^p&=&\int_\O|x|^{-p}| \psi u+(1-\psi)u|^p\\
&\leq&(1+\e)\int_\O|x|^{-p} |\psi u|^p+c(\e)\int_\O|x|^{-p}
(1-\psi)^p|u|^p \\
&\leq &(1+\e)\int_\O|x|^{-p} |\psi u|^p+ c(\e)\int_\O |u|^p.
\end{eqnarray*}
Now  by \eqref{eq:loc-H}
$$
\left(\mu_{0,p}(H)-\e\right)\int_\O|x|^{-p} |\psi u|^p\leq\int_\O|\n
(\psi u)|^p
$$
and hence
 \be\label{eq:ml-ep}
\left(\mu_{0,p}(H)-\e\right)\int_\O|x|^{-p} |u|^p\leq (1+\e)
\int_\O|\n (\psi u)|^p + c(\e)\int_\O |u|^p.
 \ee
Since $|\n (\psi u)|^p \leq \left(\psi |\n u|+ |u||\n\psi|
\right)^p$ we deduce that
$$
|\n (\psi u)|^p \leq(1+\e)\psi^p |\n u|^p+ c |u|^p|\n\psi|^p\leq
(1+\e) |\n u|^p+ c |u|^p .
$$
Using \eqref{eq:ml-ep}, we conclude that
 \be\label{eq:supgmp}
\left(\mu_{0,p}(H)-\e\right)\int_\O|x|^{-p} |u|^p\leq (1+\e)^2
\int_\O |\n u|^p+ c(\e) \int_\O|u|^p.
\ee
This implies that $\mu_{0,p}(H)\leq \sup_{\l\in
\R}\mu_{\l,p}(\O)$ and the claim follows.\\
%
%
\textbf{Step 2}: We claim that
$\sup_{\l\in\R}\mu_{\l,p}(\O)\leq\mu_{0,p}(H)$.\\
Denote by $\nu$ the unit interior normal of $\de\O$. For $\d\geq0$
we consider the cone
$$
 \C^{\d}_+:=\left\{
x\in\R^{N}~|~x\cdot\nu>\delta|x|~\right\}
$$
and put $\Sig_\d=\C^\d_+\cap \S^{N-1}$.
For every $\eta>0$, let $V\in C^\infty_c(\Sig_0)$ such that
$$
\frac{\displaystyle\int_{\Sig_0}\left(\left(\frac{N-p}{p}
\right)^2|V|^2+|\n_\s V|^2 \right)^{\frac{p}{2}}d\s}{
\displaystyle\int_{\Sig_0}|V|^pd\s}\leq \mu_{0,p}(H)+\eta.
$$%
On the other hand,  there exists $\d>0$ small such that
$supp\,V\subset \Sig_\d$. From this we conclude that
\begin{equation}\label{eq:muclmh}
 \mu_{0,p}(H)\leq \mu_{0,p}( \C^{\d}_+)\leq \mu_{0,p}(H)+\eta.
\end{equation}
Since $\de\O$ is smooth at $0$,  for every  $\d>0$, there exists
$r_\d>0$ such that
 $\C^\d_+\cap B_{r}(0)\subset
\O $  for all $r\in(0,r_\d)$. Clearly by scale invariance, $
\mu_{0,p}(\C^\d_+\cap B_{r}(0))= \mu_{0,p}(\C^\delta_+)$. For
$\e>0$, we let $\phi\in W^{1,p}_{0}( \C^\d_+\cap B_{r}(0))$ such
that
$$
\frac{\displaystyle\int_{\C^\d_+\cap B_{r}(0)}|\nabla \phi|^p~dx}
{\displaystyle\int_{\C^\d_+\cap B_{r}(0)}|x|^{-p}|\phi|^p~dx} \leq
\mu_{0,p}(\C^\delta_+)+\e.
$$
From this we deduce that
\begin{eqnarray*}
\mu_{\l,p}(\O)&\leq& \frac{\displaystyle\int_{\C^\d_+\cap
B_{r}(0)}|\nabla \phi|^p~dx-\l \int_{\C^\d_+\cap
B_{r}(0)}|\phi|^{p}~dx} {\displaystyle\int_{\C^\d_+\cap
B_{r}(0)}|x|^{-p}|\phi|^p~dx}\\
&\leq&  \mu_0(\C^\delta_+) +\e+|\l| ~ \frac{\displaystyle
\int_{\C^\d_+\cap B_{r}(0)}|\phi|^p~dx}
{\displaystyle\int_{\C^\d_+\cap B_{r}(0)}|x|^{-p}|\phi|^p~dx}.
\end{eqnarray*}
Since $\displaystyle\int_{\C^\d_+\cap B_{r}(0)}|x|^{-p}|\phi|^p~dx
\geq r^{-p} \int_{\C^\d_+\cap B_{r}(0)}|\phi|^p~dx $, we get
$$
\mu_{\l,p}(\O)\leq   \mu_{0,p}(\C^\delta_+)+\e+ r^p|\l|.
$$
The claim  follows immediately by \eqref{eq:muclmh}. Therefore
\eqref{eq:supmlmp} is proved.\\
 Finally as  the map $\l\mapsto
\mu_{\l,p}(\O)$ is non increasing while
$\mu_{\l_1,p}(\O)=0<\mu_{0,p}(H)$, we can set
$$
 \l^*(p,\O):=\inf\{\l\in\R\,:\, \mu_{\l,p}(\O)<\mu_{0,p}(H) \}
$$
so that $  \l^*(p,\O)<\mu_{0,p}(H)$ for any $\l>\l^*(p,\O)$. \QED
\begin{Remark}
 Observe that the proof of Lemma \ref{lem:supmp} highlights
that
$$
\lim_{r\to0}\mu_{0,p}(\O\cap
B_r(0))=\mu_{0,p}(H)=\lim_{\l\to-\infty}\mu_{\l,p}(\O).
$$
\end{Remark}

\bigskip
\noindent {\bf Proof of Theorem \ref{th:attained} }\\
Let $\l>\l^*(p,\O)$ so that $\mu_{\l,p}(\O)<\mu_{0,p}(H)$. We define
the mappings $F,\,G:W^{1,p}_{0}(\O)\to\R$ by
$$
F(u)=\int_\O|\n u|^p-\l\int_\O| u|^p
$$
and
$$
G(u)=\int_\O|x|^{p}| u|^p.
$$
By Ekeland variational principal, there is a minimizing sequence
$u_n\in W^{1,p}_{0}(\O)$ normalized so that
$$
G(u_n)=1,\quad \forall n\in\N
$$
and with the properties that
$$
F(u_n)\to \mu_{\l,p}(\O),
$$
\be\label{eq:almsol}
J(u_n)= F'(u_n)- \mu_{\l,p}(\O)G'(u_n)\to 0\textrm{ in }
(W^{1,p}_{0}(\O))'.
\ee
Up to a subsequence, we can  assume that  there exists $u\in
W^{1,p}_{0}(\O)$ such that
\be\label{eq:nuntonu}
\n u_n \rightharpoonup \n u \textrm{ in } L^p(\O),
 \ee
  $u_n\to u$ in $L^p(\O) $ and $ u_n\to u$ a.e. in $\O$.
   Moreover by \eqref{eq:supgmp}, we may
 assume that $|x|^{-1}u_n \rightharpoonup |x|^{-1}u$ in $L^p(\O) $. We set
 $\th_n=u_n-u$ and
$$
T(s)=\left\{ \begin{array}{ll}
s & \textrm{ if } |s|\leq 1\\
\frac{s}{|s|} & \textrm{ if } |s|> 1.
             \end{array}
       \right.
$$
It follows that for every $r\geq1$ \be\label{eq:Tzr}
\int_{\O}|T(\th_n)|^r\to 0. \ee Moreover notice that
\begin{eqnarray*}
\int_{\O}\left( |\n u_n|^{p-2}\n u_n-|\n u|^{p-2}\n u\right)\cdot \n
T(\th_n)=\la J(u_n) , T(\th_n) \ra+ \mu_{\l,p}(\O)
\int_\O|x|^{-p}|u|^{p-2}_nu_nT(\th_n)\\
+\,\l \int_\O |u|^{p-2}_n u_nT(\th_n)- \int_\O |\n u|^{p-2} \n
u\cdot \n T(\th_n).\quad
\end{eqnarray*}
Therefore by \eqref{eq:almsol}, \eqref{eq:nuntonu} and
\eqref{eq:Tzr} we infer that
$$
\int_{\O}\left( |\n u_n|^{p-2}\n u_n-|\n u|^{p-2}\n u\right)\cdot \n
T(\th_n)\to 0.
$$
Consequently by \cite{dVW}-Theorem 1.1, \be\label{eq:BLnu}
\lim_{n\to\infty}\left(\int_\O |\n u_n|^p- \int_\O |\n
\th_n|^p\right)=\int_\O |\n u|^p. \ee
By Brezis-Lieb Lemma \cite{BL} \be\label{eq:BLu} 1-
\lim_{n\to\infty}\int_\O |x|^{-p}|\th_n|^p=\int_\O |x|^{-p}| u|^p.
\ee
Fix $\e>0$ small. By  \eqref{eq:supgmp} and Rellich, there exists
$\l_\e$ such that
$$
(\mu_{0,p}(H)-\e)\int_\O|x|^{-p} |\th_n|^p\leq  \int_\O |\n
\th_n|^p- \l_\e \int_\O|\th_n|^p=  \int_\O |\n \th_n|^p+o(1).
$$
Using this together with \eqref{eq:BLnu} and \eqref{eq:BLu} we get
\begin{eqnarray*}
 \mu_{\l,p}(\O)\int_\O |x|^{-p}|u|^p&\leq& \int_\O |\n u|^p-
\l \int_\O|u|^p\leq \int_\O |\n u_n|^p-\int_\O |\n \th_n|^p-
\l \int_\O|u_n|^p+o(1)\\
&\leq &F(u_n)-\left( \mu_{0,p}(H)-\e\right)\int_\O|x|^{-p} |\th_n|^p+o(1)\\
&\leq & \mu_{\l,p}(\O) -\left( \mu_{0,p}(H)-\e\right)\left(1- \int_\O|x|^{-p} |u|^p\right)+o(1)\\
&\leq & \mu_{\l,p}(\O)- \mu_{0,p}(H)+\e +  \left(
\mu_{0,p}(H)-\e\right)\int_\O|x|^{-p} |u|^p +o(1).
\end{eqnarray*}
Send  $n\to\infty$ and then $\e\to0$ to get
$$
\left(\mu_{\l,p}(\O)- \mu_{0,p}(H) \right)\int_\O |x|^{-p}|u|^p\leq
\mu_{\l,p}(\O)- \mu_{0,p}(H).
$$
Hence  $\int_\O |x|^{-p}|u|^p\geq1$ because $\mu_{\l,p}(\O)-
\mu_{0,p}(H)<0$ and the proof is complete. \QED
As a consequence of the existence theorem, we have

\begin{Corollary}\label{cor:mrnssmmp}
 Let $\O$ be a smooth bounded domain  of $\R^N$, $N\geq2$, with $0\in\partial\Omega$.
 Then
 $$
\mu_{0,p}\left(\R^N\setminus\{0\}\right)=\left|
\frac{N-p}{p}\right|^p<\mu_{0,p}(\O)\le \mu_{0,p}(H)~\!.
 $$
\end{Corollary}

\proof By \eqref{eq:supmlmp}   $0<\mu_{0,p}(\O)\le \mu_{0,p}(H)$. If
the strict inequality holds, then there exists a positive minimizer
$u\in W^{1,p}_0(\O)$ for $\mu_{0,p }(\O)$ by Theorem
\ref{th:attained}. But then
$\mu_{0,p}\left(\R^{N}\setminus\{0\}\right)<\mu_{0,p}(\O)$, because
otherwise a null extension of $u$ outside $\O$ would achieve the
Hardy constant in $\R^N\setminus\{0\}$ which is not possible. \QED
As mentioned earlier, we shall  show that there are smooth bounded
domains in $\R^N$ such that $\l^*(p,\O)\in[-\infty,0)$. These
domains might be taken to be convex or even flat at $0$. For that we
let $\nu\in\S^{N-1}$ and $\d,r,R>0$. We consider the sector
\begin{equation}\label{eq:sect}
 \calC^{\d}_{r,R}:=\left\{
x\in\Rnk~|~x\cdot\nu>-\delta|x|~\,, r<|x|<R\right\}.
\end{equation}
\begin{Proposition}\label{prop:uppls}
Let  $N\geq2$ and $p>1$. Then for all $\d\in(0,1)$, there exist
$r,R>0$ such that
if a domain $\O$ contains $\calC^{\d}_{r,R}$ then
$\mu_{0,p}(\O)<\mu_{0,p}(H)$.
\end{Proposition}
\proof
Consider the cone
$$
 \calC^{\d}:=\left\{
x\in\Rnk~|~x\cdot\nu>-\delta|x|~\right\}
$$
 Notice that by Harnack inequality
 $\mu(\calC^{\d})<\mu(\calC^{\d'})$ for any $0\leq \d'<\d<1$. Thus
for any $\d\in(0,1)$, we can find $u\in C^\infty_c(\calC^\d)$ such
that
 $$
 \frac {\displaystyle\int_{\calC^\d}|\n
 u|^p}{\displaystyle\int_{\calC^\d}|x|^{-p}|u|^p}<\mu_{0,p}(H).
 $$
 Hence we choose $r,R>0$ so that
 ${supp}~u\subset\calC^\d_{r,R}$.
 \QED
%
%
By Corollary \ref{cor:mrnssmmp}, starting from \textrm{exterior
domains}, one can also build various example of (possibly annular)
domains for which $\l^*(p,\O)<0$. The following argument is taken in
[Ghoussoub-Kang~\cite{GK} Proposition~2.4]. If $U\subset\R^N$,
$N\geq2,$ is a smooth exterior domain (the complement of a smooth
bounded domain) with $0\in\de U$ then by scale invariance
$\mu_{0,p}(U)=\mu_{0,p}(\R^{N}\setminus\{0\})$. We let
 $B_r(0)$ a ball of radius $r$ centered at the 0 and define $\O_r:= B_r(0)\cap U$
  then clearly the
map $r\mapsto \mu(\O_r)$ is decreasing with
\begin{equation}\label{eq:lssmlmp}
\mu_{0,p}\left(\R^{N}\setminus\{0\}\right)=\inf_{r>0}
\mu_{0,p}(\O_r)\quad\textrm{ and }\quad \mu_{0,p}(H)=\sup_{r>0}
\mu_{0,p}(\O_r).
\end{equation}
We have the following result for which the proof is similar  to the
one given in \cite{GK} by Corollary \ref{cor:mrnssmmp} and Harnack
inequality.
\begin{Proposition}\label{prop:uppls1}
There exists $r_0>0$ such that the mapping $r\mapsto
\mu_{0,p}(\O_r)$ is left-continuous and strictly decreasing on
$(r_0,+\infty)$. In particular
$$\mu_{0,p}(\R^{N}\setminus\{0\})<\mu_{0,p}(\O_r)<\mu_{0,p}(H),\quad\forall r\in (r_0,+\infty).$$
\end{Proposition}
\subsection{Remainder term}\label{ss:rt}
We know that for domains $\O$ contained in a half-ball
$\l^*(p,\O)\geq0$. Our aim in this section is to obtain positive
lower bound for $\l^*(p,\O)$ by providing a remainder term for
Hard's inequality in these domains.
In \cite{GGM}, Gazzola-Grunau-Mitidieri proved  the following
improved Hardy inequality for $1<p<N$:
\begin{equation}
\label{eq:BV} \int_\Omega|\nabla u|^p-
\mu_{0,p}\left(\R^{N}\setminus\{0\}\right)\int_\Omega|x|^{-p}|u|^p\ge
C(N,p)\left(
\frac{\o_N}{|\Omega|}\right)^\frac{p}{N}\int_\Omega|u|^p~\!,
\end{equation}
that holds for any bonded domain $\O$ of $\R^N$ and  $u\in
W^{1,p}_0(\O)$. Here the constant $C(N,p)>0$ is explicitly given
while  $C(N,2)$ is the first Dirichlet eigenvalue of $-\D$ of the
unit disc in $\R^2$.
\\
We shall show that such  type of inequality holds in the case where
the singularity is placed at the boundary of the domain. To this
end, we will use the function $v(x):=|x|^{\frac{p-N}{p}}
V\left(\frac{x}{|x|}\right)$ defined in \eqref{eq:Hardeigf} to
"reduce the dimension".\\
 Throughout this section, we assume that
$N\geq2$
 since the case $N=1$ was already proved by Tibodolm \cite{Tib} Theorem 1.1.
  Indeed, he showed that \\
 $$
\int_0^1|u'(r)|^pdr-\mu_{0,p}(H)\int_0^1r^{-p}|u(r)|^pdr\geq(p-1)2^{p}
\int_0^1|u(r)|^pdr,\quad\forall u\in W^{1,p}_0(0,1).
 $$
 We start with conic domains
$$
\C_{\Sigma}=\{x= r\sigma\in\R^N~|~r\in(0,1),\, \sigma\in \Sigma~\},
$$
where $\Sigma$ is a domain properly contained in $\S^{N-1}$ and
having a Lipschitz boundary. We will denote by $V$ the positive
minimizer of \eqref{eq:mpH} in $\Sig$ while
$v(x):=|x|^{\frac{p-N}{p}} V\left(\frac{x}{|x|}\right)$ satisfies
\eqref{eq:vweaksol} in the infinite cone $\{x=
r\sigma\in\R^N~|~r\in(0,+\infty),\, \sigma\in \Sigma~\}$. Finally we
remember that by Harnack inequality $\frac{1}{v}\in
L^\infty_{loc}(\calC_\Sig)$.

Recall the following inequalities (see \cite{Lind} Lemma 4.2) which
will be useful in the remaining of the paper. Let $p\in [2,\infty)$
then for any $a,b\in\R^N$
\be\label{eq:apbp}
|a+b|^p\geq |a|^p+\frac{1}{2^{p-1}-1}|b|^p+p|a|^{p-2}a\cdot b.
\ee
If  $p\in (1,2)$ then for any $a,b\in\R^N$
\be\label{eq:apbp2}
|a+b|^p\geq
|a|^p+c(p)\frac{|b|^2}{\left(|a|+|b|\right)^{2-p}}+p|a|^{p-2}a\cdot
b.
\ee

 We first make the following observation.
\begin{Lemma}
Let $u\in C^\infty_c(\calC_{\Sig})$, $u\geq0$. Set $\psi=\frac{u}{v}$ then \\
If $p\geq2$
\be\label{eq:impg2}
\int_{\calC_{\Sig}}|\n
u|^p-\mu_{0,p}(\calC_{\Sig})\int_{\calC_{\Sig}}|x|^{-p}|u|^p\geq\frac{1}{2^{p-1}-1}\int_{\calC_{\Sig}}|v\n\psi|^p,
\ee
If $1<p<2$
\be \label{eq:impl2}
\int_{\calC_{\Sig}}|\n
u|^p-\mu_{0,p}(\calC_{\Sig})\int_{\calC_{\Sig}}|x|^{-p}|u|^p\geq
c(p)\int_{\calC_{\Sig}}\frac{|v\n\psi|^2}{\left(|v\n\psi|+|\psi\n v|
  \right)^{2-p}},
\ee
\end{Lemma}
\proof
We prove only the case $p\geq2$ as the  case $p\in(1,2)$ goes
similarly. Notice that $\n u=v\n \psi+\psi\n v$ then  we use the
inequality \eqref{eq:apbp}  with $a=v\n \psi$ and $b=\psi\n v$ to
get
$$
\int_{\calC_{\Sig}}|\n u|^p\geq \int_{\calC_{\Sig}}|\psi \n
v|^p+p\int_{\calC_{\Sig}}|\psi\n v|^{p-2}\psi\n
v\cdot(v\n\psi)+\frac{1}{2^{p-1}-1}\int_{\calC_{\Sig}}|v\n\psi|^p.
$$
It is plain that
$$
p|\psi\n v|^{p-2}\psi\n v\cdot(v\n\psi)=|\n v|^{p-2}\n
v\cdot(v\n\psi^p) =|\n v|^{p-2}\n v\cdot\n(v\psi^p) - |\psi\n
v|^{p}.
$$
Inserting this in the first inequality and using \eqref{eq:vweaksol}
we deduce that
\begin{eqnarray*}
 \int_{\CS}|\n u|^{p}&\geq&
\frac{1}{2^{p-1}-1}\int_{\CS}|v\n\psi|^p+\int_{\CS} |\n
v|^{p-2}\n v\cdot\n(v\psi^p)\\
&\geq& \frac{1}{2^{p-1}-1}\int_{\CS}|v\n\psi|^p+\mu_{0,p}(\CS)
\int_{\CS}|x|^{-p}u^p.
\end{eqnarray*}
 \QED
The improvement in the case $p\geq2$ is an immediate consequence of
the above lemma.
\begin{Lemma}\label{lem:rempgeq2} For all
$p\geq2$
 $$
 \int_{\calC_{\Sig}}|\n
u|^p-\mu_{0,p}(\calC_{\Sig})\int_{\calC_{\Sig}}|x|^{-p}|u|^p\geq\frac{\L_p}{2^{p-1}-1}\int_{\CS}|u|^p,\quad\forall
u\in C^\infty_c(\CS),
$$
where $ \L_p:=\inf_{f\in
C^1_c(0,1)}\frac{\int_0^1r^{p-1}|f'|^pdr}{\int_0^1r^{p-1}|f|^pdr}$.
\end{Lemma}
\proof Since $|\n|u||\leq|\n u|$, we may assume that $u\geq0$.
 We only need to estimate the right hand side in
\eqref{eq:impg2}. We use polar coordinates $x\mapsto
(|x|,\frac{x}{|x|})=(r,\s)$ and denote by $\de_r$ the radial
direction. Then using \eqref{eq:apbp},
\begin{eqnarray*}
\int_{\calC_{\Sig}}|v\n\psi|^p&=&\int_{\Sig}\int_0^1r^{p-1}V^p|\psi_r\de_r+\n_\s\psi|^p\\
&\geq& \int_{\Sig}V^p \int_0^1r^{p-1}|\psi_r|^p\geq\L_p
\int_{\Sig}V^p\int_0^1r^{p-1} |\psi|^p\\
&\geq&\L_p \int_{\Sig}\int_0^1u^pr^{N-1}=\L_p\int_{\CS}|u|^p.
\end{eqnarray*}
The lemma readily follows from \eqref{eq:impg2}.
 \QED

 It is easy to see that by integration by parts $\L_p\geq1$
while for integer $p\in\N$ then $\L_p$ corresponds to the first
Dirichlet eigenvalue of $-\D$ in the unit ball of $\R^p$.\\
We now turn to the case $p\in (1,2)$ which carries more
difficulties. We shall need the following intermediate result.
\begin{Lemma}\label{lem:upbrr2}
Let $p\in (1,2)$ and $u\in C^\infty_c(\CS)$, $u\geq0$. Setting
$\psi=\frac{u}{v}$ then there exists a constant $c=c(p,\Sig)>0$ such
that
$$
c\int_{\CS}r|\psi\n v|^p\leq \int_{\CS}r^{(2-p)/2}|v\n\psi|^p.
$$
\end{Lemma}
\proof Let $\ti{\psi}:={r^{\frac{1}{p}}}\psi$ and use  $\ti{\psi}^p
v$ as a test function in the weak equation \eqref{eq:vweaksol}. Then
by H\"{o}lder
\begin{eqnarray*}
\int_{\CS}|\ti{\psi}\n v|^p&\leq&
\mu_{0,p}(\CS)\int_{\CS}r^{-p}v^p\ti{\psi}^p+p\int_{\CS}|\ti{\psi}\n
v|^{p-1}|v\n\ti{\psi}|\\
&\leq&\mu_{0,p}(\CS)\int_{\CS}r^{-p}v^p\ti{\psi}^p+p\left(\int_{\CS}|\ti{\psi}\n
v|^p  \right)^{\frac{p-1}{p}}\left(\int_{\CS}|v\n\ti{\psi}|^p  \right)^{\frac{1}{p}}.\\
\end{eqnarray*}
Therefore by Young's inequality, for $\e>0$ small there exists a
constant $C_\e>0$ depending on $p$ and $\Sig$ such that
$$
(1-\e c(p))\int_{\CS}|\ti{\psi}\n v|^p\leq C_\e
\int_{\CS}r^{-p}v^p\ti{\psi}^p +C_\e\int_{\CS}|v\n\ti{\psi}|^p .
$$
Recall that  $\ti{\psi}=r^{\frac{1}{p}}\psi$. Then since
$$
|\n\ti{\psi}|^p\leq c(p)\left(r^{1-p}\psi^p+r|\n\psi|^p \right),
$$
we conclude that there exists a constant $c=c(p,\Sig)$ such that
\be\label{eq:befH2}
c\int_{\CS}r|{\psi}\n v|^p\leq \int_{\CS}r^{1-p}v^p{\psi}^p
+\int_{\CS}r^{(2-p)/p}|v\n{\psi}|^p,
 \ee
 we have used the fact that $r\leq r^{(2-p)/p}$ for all $r\in(0,1)$. To
 estimate the first term in the right hand side in \eqref{eq:befH2}
 we will use the 2-dimensional Hardy inequality. Through the polar
 coordinates $x\mapsto (r,\s)$
\begin{eqnarray*}
\int_{\CS}r^{1-p}v^p{\psi}^p&=&\int_{\Sig}V^p\int_0^1r^{p-1}\left(\frac{\psi}{r}\right)^pr\\
&\leq&\int_{\Sig}V^p\int_0^1\left(\frac{\psi}{r}\right)^pr\\
&\leq& \left| \frac{p}{p-2}\right|^{-p}
\int_{\Sig}V^p\int_0^1|\psi_r|^pr\\
&=&\left| \frac{p}{p-2}\right|^{-p} \int_0^1\int_{\Sig}V^p
v^{-p}|v\n\psi|^pr=\left|
\frac{p}{p-2}\right|^{-p}\int_0^1\int_{\Sig} r^{N-p+1}|v\n\psi|^p.
\end{eqnarray*}
To conclude, we notice that
$r^{N-p+1}=r^{N-\frac{p}{2}}r^{(2-p)/p}\leq r^{N-1}r^{(2-p)/p}$ as
$p\in (1,2)$ so that
$$
\int_{\CS}r^{1-p}v^p{\psi}^p\leq \left|
\frac{p}{p-2}\right|^{-p}\int_{\CS}r^{(2-p)/p}|v\n\psi|^p.
$$
Inserting this in \eqref{eq:befH2} the lemma  follows immediately.
 \QED
We are now in position to prove the improved Hardy inequality for
$p\in(1,2)$.
\begin{Lemma}\label{lem:remple2}
Let $p\in(1,2)$.
Then there exists a  constant $c=c(p,\Sig)>0$ such that
 $$
 \int_{\calC_{\Sig}}|\n
u|^p-\mu_{0,p}(\calC_{\Sig})\int_{\calC_{\Sig}}|x|^{-p}|u|^p\geq
c\int_{\CS}|u|^p,\quad\forall u\in C^\infty_c(\CS).
$$
\end{Lemma}
\proof Here also we may assume that $u\geq0$. We need to estimate
the right hand side of \eqref{eq:impl2}. Let  $r=|x|$ then by
H\"{o}lder and Lemma \ref{lem:upbrr2}, we  have
\begin{eqnarray*}
\int_{\CS}r^{\frac{2-p}{2}}|v\n\psi|^p&=&\int_{\CS}\frac{|v\n\psi|^2}{\left(|v\n\psi|+|\psi\n
v| \right)^{(2-p)p/2}}   r^{\frac{2-p}{2}}\left(|v\n\psi|+|\psi\n v|
\right)^{(2-p)p/2}\\
&\leq&\left(   \int_{\CS}\frac{|v\n\psi|^2}{\left(|v\n\psi|+|\psi\n
v| \right)^{2-p}} \right)^{p/2}   \left(\int_{\CS}r \left|
|v\n\psi|+|\psi\n v|  \right|^p \right)^{ (2-p)/2}\\
&\leq& \left(   \int_{\CS}\frac{|v\n\psi|^2}{\left(|v\n\psi|+|\psi\n
v| \right)^{2-p}} \right)^{p/2} \\
&\,&\times \left(2^{p-1}\int_{\CS}r
|v\n\psi|^p+2^{p-1}\int_{\CS}r|\psi\n v|^p \right)^{ (2-p)/2}\\
&\leq& c \left( \int_{\CS}\frac{|v\n\psi|^2}{\left(|v\n\psi|+|\psi\n
v| \right)^{2-p}} \right)^{p/2} \left(\int_{\CS}r^{\frac{2-p}{2}}
|v\n\psi|^p \right)^{ (2-p)/2},
\end{eqnarray*}
where $c$ a positive constant depending only on $p$ and $\Sig$ and
we have used once more the fact that $r\leq r^{(2-p)/p}$ for all
$r\in(0,1)$. Consequently by \eqref{eq:impl2}, we deduce that
\be\label{eq:almosimpH}
\int_{\calC_{\Sig}}|\n
u|^p-\mu_{0,p}(\calC_{\Sig})\int_{\calC_{\Sig}}|x|^{-p}|u|^p\geq
c\int_{\CS}r^{\frac{2-p}{2}}|v\n\psi|^p.
\ee
To  proceed  we  estimate
\begin{eqnarray*}
\int_{\Sig}\int_0^1u^pr^{N-1}&=&\int_{\Sig}V^p\int_0^1r^{p-1}
|\psi|^p\leq c(p) \int_{\Sig}V^p\int_0^1r|\psi_r|^p\\
&\leq &c(p)\int_{\Sig}V^p\int_0^1r^{\frac{p}{2}}|\psi_r|^p\\
&\leq&c(p)\int_{\CS}r^{\frac{2-p}{2}}|v\n\psi|^p.
\end{eqnarray*}
The first inequality comes from  the 2-dimensional embedding
$W^{1,p}_{0}\subset L^{\frac{2p}{2-p}}\subset L^{\frac{p}{3-p}}$,
one can see [\cite{GGM} page 2155] for the proof.
Putting this in \eqref{eq:almosimpH} we conclude  that there exists
a positive constant $c=c(p,\Sig)$ such that
$$
\int_{\calC_{\Sig}}|\n
u|^p-\mu_{0,p}(\calC_{\Sig})\int_{\calC_{\Sig}}|x|^{-p}|u|^p\geq
c\int_{\CS}|u|^p
$$
which was the purpose of the lemma.
 \QED
The main result in this section is contained in the next theorem.

\begin{Theorem}\label{th:bd-ls}
Let $\O$ be a domain in  $\R^N$ with $0\in\de\Omega$. If $\O$ is
contained in a half-ball centered at 0 then there exists a constant
$c(N,p)>0$ such that
$$
\int_{\O}|\nabla u|^p-\mu_{0,p}(H)\int_{\O} |x|^{-p}|u|^p\geq
\frac{c(N,p)}{\textrm{diam}(\O)^p}\int_{\O} |u|^p\qquad \forall u\in
W^{1,p}_{0}(\O).
$$
\end{Theorem}

\proof Let $R=\textrm{diam}(\O)$ be the diameter of $\Omega$. Then
$\Omega$ is contained  in a half  ball  $B^+_R$ of radius $R$
centered at the origin. From Lemma \ref{lem:rempgeq2} and Lemma
\ref{lem:remple2} we infer that
 $$
\int_{B^+_R}|\nabla u|^p-\mu_{0,p}(H)\int_{B^+_R} |x|^{-p}|u|^p\geq
\frac{c(N,p)}{R^p}\int_{B^+_R}| u|^p\quad \forall u\in
C^{\infty}_c(\O)
$$
by homogeneity. The theorem readily follows by density. \QED

We do not know whether $\textrm{diam}(\O)$ might be replaced with
$\o_N|\O|^\frac{1}{N}$ as in \cite{GGM} at least when $\O$ is convex
and $p\geq2$. There might exists also  "logarithmic" improvement as
was recently obtained in \cite{FaMu1} inside  cones and  $p=2$. One
can see also the work of Barbatis-Filippas-Tertikas in \cite{BFT}
for domains containing the origin or when $|x|$ is replaced by the
distance to the boundary.

\appendix
\section{Hardy's inequality}\label{app:Hd}
We denote by ${d}$ the distance function of $\O$:
$$
{d}(x):=\inf\{|x-\s|\, :\, \s\in \de\O\}.
$$
In this section, we study the  problem of finding minima to the
following quotient
\begin{equation}
\label{eq:problem-n} \nu_{\l,p}(\O):= \inf_{u\in W^{1,p}_{0}(\O)}
~\frac{\displaystyle\int_{\O}|\nabla u|^p~dx-\l\int_{\O}|u|^p~dx}
{\displaystyle\int_{\O}{d}^{-p}|u|^p~dx}~,
\end{equation}
where $p>1$ and  $\lambda\in\R$ is a varying parameter. Existence of
extremals to this problem was studied in \cite{BM} when $p=2$ and in
\cite{MMP} with $\l=0$. It is known (see for instance \cite{MMP})
that $\nu_{0,p}(\O)\leq \textbf{c}_p$ for any smooth bounded domain
$\O$ while  for convex domain $\O$, the Hardy constant
$\nu_{0,p}(\O)$   is not achieved and
$\nu_{0,p}(\O)=\left(\frac{p-1}{p}\right)^p=:\textbf{c}_p$.
\\
The main  result in this section is contained in the following
\begin{Theorem}\label{th:attained-d}
 Let $\O$ be a smooth
bounded domain in $\R^N$  and  $p>1$, there exits $\ti{\l}(p,\O)
\in[-\infty,+\infty)$ such that
\be
\nu_{\l,p}(\O)<\left(\frac{p-1}{p}\right)^p,\quad \forall
\l>\ti{\l}(p,\O).
\ee
The infinimum in \eqref{eq:problem-n}  is attained if
$\l>\ti{\l}(p,\O)$.
\end{Theorem}


%

  %
  We start with the following result which is stronger than needed.
  It was proved in \cite{BM} for $p=2$ and in \cite{FMT} when $2\leq p<N$ as
  the authors were dealing with Hardy-Sobolev inequalities.
  \begin{Lemma}
 Let $\O$ be a smooth bounded domain in $\R^N$ and
 $p\in(1,\infty)$. Then there exists $\b=\b(p,\O)>0$ small such that
 \be\label{eq:HardyOb}
 \int_{\O_\b}|\n u|^p\geq \textbf{c}_p\int_{\O_\b}{d}^{-p}|u|^p\quad\forall u\in
 H^1_0(\O),
 \ee
 where  $
\O_\b:=\{x\in\O\,:\,{d}(x)<\b\}.
$
\end{Lemma}
\proof Since $|\n|u||\leq|\n u|$, we may assume that $u\geq0$. Let
$u\in C^\infty_c(\O)$ and put $v={d}^{\frac{1-p}{p}}u$. Using
\eqref{eq:apbp} and \eqref{eq:apbp2}, we get
\begin{equation}\label{eq:pg2d}
|\n u|^p-\textbf{c}_p{d}^{-p}|u|^p\geq c(p)d^{p-1}|\n v|^p+
\left|\frac{p-1}{p}\right|^{p-1}\n{d}\cdot\n (v^{p})\quad\textrm{ if
} p\geq2,
\end{equation}
\begin{equation}\label{eq:pl2d}
 |\n u|^p-\textbf{c}_p{d}^{-p}|u|^p\geq c(p)\frac{{d}|\n
v|^2}{\left(\textbf{c}_p^{\frac{1}{p}}|v|+{d} |\n v|\right)^{2-p}}+
\left|\frac{p-1}{p}\right|^{p-1}\n{d}\cdot\n (v^{p})\quad\textrm{ if
} p\in(1,2).
\end{equation}
By integration by parts, we have
$$
\int_{\O_{\b}}\n{d}\cdot\n (v^{p})=-\int_{\O_\b}\D {d}|v|^p+
\int_{\de\O_\b}|v|^{p}\geq
-c\int_{\O_\b}|v|^{p}+\int_{\de\O_\b}|v|^{p},
$$
for a positive constant depending only on $\O$. Multiply the
identity $\div({d}\n{d})=1+{d}\D{d}$ by $v$ in integrate by parts to
get
$$
\left(1+o(1)
\right)\int_{\O_\b}|v|^{p}=-p\int_{\O_\b}{d}|v|^{p-1}\n{d}\cdot\n
v+\int_{\de \O_\b}{d} |v|^{p}\leq c(p)\int_{\O_\b}{d}|v|^{p-1}|\n
v|+\int_{\de \O_\b}{d} |v|^{p}.
$$
By H\"{o}lder and Young's inequalities
\be\label{eq:estv}
\left(1+o(1)-c\e \right)\int_{\O_\b}|v|^{p}\leq c_\e \int_{
\O_\b}{d}^p |\n v|^{p}+\int_{\de \O_\b}{d} |v|^{p}.
\ee
\textbf{Case  $p\geq 2$}. Using \eqref{eq:estv} we infer that
$$
\left(1+o(1)-c\e \right)\int_{\O_\b}|v|^{p}\leq c_\e \b \int_{
\O_\b}{d}^{p-1} |\n v|^{p}+\b \int_{\de \O_\b} |v|^{p}.
$$
It follows from \eqref{eq:pg2d} that for $\e,{\b}>0$ small
$$
\int_{\O_\b}|\n u|^p-\textbf{c}_p\int_{\O_\b}{d}^{-p}|u|^p\geq
c\left( \int_{\O_\b}{d}^{p-1}|\n v|^p+\int_{\de \O_\b}
|v|^{p}\right)
$$
as desired.\\
 \textbf{Case  $p\in(1,2)$}. By H\"{o}lder and Young's
inequalities
\begin{eqnarray*}
\int_{\O_\b}{d}^{p}|\n
v|^p&=&\displaystyle\int_{\O_\b}{\frac{{d}^p|\n
v|^p}{\left(\textbf{c}_p^{\frac{1}{p}}|v|+{d}|\n v|
\right)^{\frac{p(2-p)}{2}}}\left(
\textbf{c}_p^{\frac{1}{p}}|v|+{d}|\n
v|\right)^{\frac{p(2-p)}{2}}}\\
&\leq&c_\e \displaystyle\int_{\O_\b}\frac{{d}^2|\n
v|^2}{\left(\textbf{c}_p^{\frac{1}{p}}|v|+{d}|\n v|
\right)^{2-p}}+\e c\int_{\O_\b}|v|^p+\e c \int_{\O_\b}{d}^p|\n v|^p
\end{eqnarray*}
and thus
$$
\left(1-c\e  \right)\int_{\O_\b}{d}^{p}|\n v|^p\leq c_\e
\displaystyle\int_{\O_\b}\frac{{d}^2|\n
v|^2}{\left(\textbf{c}_p^{\frac{1}{p}}|v|+{d}|\n v|
\right)^{2-p}}+\e c\int_{\O_\b}|v|^p.
$$
Using this in \eqref{eq:estv} we obtain
$$
\left(1+o(1)-c\e \right)\int_{\O_\b}|v|^{p}\leq c\b
\displaystyle\int_{\O_\b}\frac{{d}|\n
v|^2}{\left(\textbf{c}_p^{\frac{1}{p}}|v|+{d}|\n v|
\right)^{2-p}}+c\b\int_{\de \O_\b} |v|^{p} .
$$
By \eqref{eq:pl2d}, we conclude that for $\e,\b>0$ small
$$
\int_{\O_\b}|\n u|^p-\textbf{c}_p\int_{\O_\b}{d}^{-p}|u|^p\geq
\displaystyle c\int_{\O_\b}\frac{{d}|\n
v|^2}{\left(\textbf{c}_p^{\frac{1}{p}}|v|+{d}|\n v|
\right)^{2-p}}+c\int_{\de \O_\b} |v|^{p}.
$$
This ends the proof of the lemma.
 \QED

 \begin{Lemma}\label{lem:supnup}
  Let $\O$ be a bounded domain of class $C^2$ in $\R^N$. Then there exists
  $\ti{\l}(p,\O)\in[-\infty,+\infty)$ such that
  $$
  \nu_{\l,p}(\O)<\textbf{c}_p\quad\forall \l>\ti{\l}(p,\O).
  $$
  \end{Lemma}
\proof
The proof will be carried out in 2 steps.\\
 \textbf{Step 1}: We claim
that
$\sup_{\l\in\R}\nu_{\l,p}(\O)\geq\textbf{c}_p$.\\
For $\b>0$  we define
$$
\O_\b:=\{x\in\O\,:\,{d}(x)<\b\}.
$$
  Let $\psi\in C^\infty(\O_\b)$ with $0\leq\psi\leq 1$, $\psi\equiv
0$ in $\R^N\setminus \O_{\frac{\b}{2}}$ and $\psi\equiv 1$ in
$\O_{\frac{\b}{4}}$. For $\e>0$ small, there holds
\begin{eqnarray*}
\int_\O{d}^{-p} |u|^p&=&\int_\O{d}^{-p}| \psi u+(1-\psi)u|^p\\
&\leq&(1+\e)\int_\O{d}^{-p} |\psi u|^p+C\int_\O{d}^{-p}
(1-\psi)^p|u|^p \\
&\leq &(1+\e)\int_\O{d}^{-p} |\psi u|^p+ C\int_\O |u|^p.
\end{eqnarray*}
By \eqref{eq:HardyOb}, we infer that
$$
\textbf{c}_{p}\int_\O{d}^{-p} |\psi u|^p\leq\int_\O|\n (\psi u)|^p
$$
and hence
 \be\label{eq:ml-epd}
\textbf{c}_{p}\int_\O{d}^{-p} |u|^p\leq (1+\e) \int_\O|\n (\psi
u)|^p + C\int_\O |u|^p.
 \ee
Since $|\n (\psi u)|^p \leq \left(\psi |\n u|+ |u||\n\psi|
\right)^p$ we deduce that
$$
|\n (\psi u)|^p \leq(1+\e)\psi^p |\n u|^p+ C |u|^p|\n\psi|^p\leq
(1+\e) |\n u|^p+ C |u|^p .
$$
Using \eqref{eq:ml-epd}, we conclude that
$$
\textbf{c}_{p}\int_\O{d}^{-p} |u|^p\leq (1+\e)^2 \int_\O |\n u|^p+
C(\e,\b)\, \int_\O|u|^p.
$$
This means that $\textbf{c}_{p}\leq \sup_{\l\in \R}\nu_{\l,p}(\O)$.\\
 \textbf{Step 2}: We claim
that
$\sup_{\l\in\R}\nu_{\l,p}(\O)\leq\textbf{c}_p$.\\
%
%
Let $\b>0$ then by $\eqref{eq:Hardy-1d}$ and scale invariance we
have $\mu_{0,p}(0,\b)=\textbf{c}_p$. Hence  for $\e>0$ there exits a
function $\phi\in W^{1,p}_{0}(0,\b)$ such that
\be\label{eq:1dH} \textbf{c}_p+\e\geq \frac{ \int_0^\b|\phi'|^p \,
ds }{ \int_0^\b s^{-p}\phi^p\, ds}.
\ee
Letting $u(x)=\phi({d}(x))$, there exists a positive constant $C$
depending only on $\O$ such that
$$
\int_{\O_\b}|\n u|^p=\int_0^\b\int_{\de\O_s}|\phi'(s)|^p\, d\s_s
 \leq \left( 1+C\b\right)|\de\O|\int_0^\b|\phi'(s)|^p\,ds.
$$
Furthermore
$$
\int_{\O_\b}{d}^{-p}|
u|^p=\int_0^\b\int_{\de\O_s}s^{-p}|\phi(s)|^p\, d\s_s \geq \left(
1-C\b\right)|\de\O|\int_0^\b|\phi(s)|^p\, ds.
$$
By \eqref{eq:1dH} we conclude that

$$
\nu_{\l,p}(\O)\leq \frac{\displaystyle\int_{\O_\b}|\nabla u|^p~dx-\l
\int_{\O_\b}u^{p}~dx} {\displaystyle\int_{\O_\b}{d}^{-p}|u|^p~dx}
\leq ({\textbf{c}_p}+\e)\frac{1+C\b}{1-C\b}+|\l| \frac{\displaystyle
\int_{\O_\b}|u|^{p}~dx} {\displaystyle\int_{\O_\b}{d}^{-p}|u|^p~dx}.
$$
Since $\int_{\O}{d}^{-p}|u|^2~dx \geq \b^{-p} \int_{\O_\b}|u|^p~dx
$, we get
$$
\nu_{\l,p}(\O)\leq  ({\textbf{c}_p}+\e)\frac{1+C\b}{1-C\b} +
\b^p|\l|,
$$
sending $\b$ to $0$ we get the desired result.
\QED
Clearly the proof of Theorem~\ref{th:attained-d} goes similarly as
the one
 of Theorem~\ref{th:attained} and we skip it.\\
It was shown in \cite{BM} that $\ti{\l}(2,\O)\in\R$ and that
$\nu_{\l,p}(\O) $ is not achieved for any $\l\geq \ti{\l}(2,\O)$. On
the other hand by \cite{Davis}, there are domains for which
$\ti{\l}(2,\O)<0$,
see also \cite{MMP}.\\
 We point out that if $\O$ is convex then by \cite{Tib} there
exists a constant $a(N,p)>0$ (explicitly given) such that
$$
\ti{\l}(p,\O)\geq\frac{a(N,p)}{|\O|^{\frac{p}{N}}}.
$$
 We finish this section by showing that there are
smooth bounded domains in $\R^N$ such that
$\ti{\l}(p,\O)\in[-\infty,0)$. We let $U\subset\R^N$, $N\geq2$ with
$0\in\de U$ be an exterior domain and set $\O_r=B_r(0)\cap U$.
\begin{Proposition}
Assume that $p>\frac{N+1}{2}$ then there exists $r>0$ such that $
\nu_{0,p}(\O_r)<\left(\frac{p-1}{p}\right)^p.$
\end{Proposition}
\proof
Clearly  $\mu_{0,p}(\R^N\setminus\{0\})= \left|
\frac{N-p}{p}\right|^p< \left(\frac{p-1}{p}\right)^p $ provided
$p>\frac{N+1}{2}$.  Let $\e>0$  such that
$\left(\frac{p-1}{p}\right)^p>\left| \frac{N-p}{p}\right|^p+\e$ so
by \eqref{eq:lssmlmp}, there exits ${r}>0$ such that
$$
\mu_{0,p}(\O_r)< \left|
\frac{N-p}{p}\right|^p+\e<\left(\frac{p-1}{p}\right)^p.
$$
The conclusion readily follows  since $\nu_{0,p}(\O_r)\leq
\mu_{0,p}(\O_r)$ because $0\in\de\O_r$.
\QED

%
%
%
%
%
%
%
%

%
%

%

\begin{thebibliography}
\footnotesize
%
\bibitem{BFT}   Barbatis G., Filippas S.,  Tertikas A., A unified approach to
improved $L^p$ Hardy inequalities with best constants . Trans. Amer.
Math. Soc., 356, (2004), 2169-2196.
%
\bibitem{BM} Brezis H. and Marcus M., Hardy's inequalities revisited.
 Dedicated to Ennio De Giorgi.
  Ann. Scuola Norm. Sup. Pisa Cl. Sci. (4)  25  (1997),  no. 1-2, 217-237.
%
\bibitem{BV} Brezis H. and V\`{a}zquez J. L., Blow-up solutions of some nonlinear elliptic problems.
 Rev. Mat. Univ. Complut. Madrid  10  (1997),  no. 2, 443-469.

%
%
\bibitem{BL} Brezis H. and  Lieb E., A relation between pointwise convergence of functions
 and convergence of functionals,
 Proc. Amer. Math. Soc. 88 (1983), 486-490.
%
\bibitem{CaMuPRSE}
{  Caldiroli P.,  Musina R.}, On a class of 2-dimensional singular
elliptic problems. Proc. Roy. Soc. Edinburgh Sect. A 131 (2001),
479-497.
%

\bibitem{CaMuUMI} Caldiroli P., Musina R., Stationary states for a
two-dimensional~ singular Schr\"odinger equation. Boll. Unione Mat.
Ital. Sez. B Artic. Ric. Mat.  (8) 4-B (2001), 609-633.

%
\bibitem{dVW} de Valeriola S. and Willem M., On Some Quasilinear Critical Problems,
 Advanced Nonlinear Studies 9 (2009), 825-836.
%
\bibitem{Davis} Davies E. B., The Hardy constant, Quart. J. Math. Oxford (2) 46
(1995), 417-431.
%
\bibitem{MMF-PJM} Fall M. M., Area-minimizing regions with small volume
in Riemannian manifolds with boundary. Pacific J. Math. 244 (2010),
no. 2, 235-260.
%
\bibitem{FaMu} Fall M. M., Musina R., Hardy-Poincar\'e inequalities with boundary singularities.
Pr\'epublication D\'epartement de Math\'ematique Universit\'e
Catholique de Louvain-La-Neuve 364 (2010), {\tt
http://www.uclouvain.be/38324.html}.
%
\bibitem{FaMu1} Fall M. M., Musina R., Sharp nonexistence results for a linear elliptic inequality
involving Hardy and Leray potentials. Prerint SISSA (2010). Ref.
31/2010/M.
%
\bibitem{FMT}  Filippas S., Maz'ya V. and Tertikas A., Critical Hardy$-$Sobolev inequalities.
Journal de Math\'ematiques Pures et Appliqu\'es Volume 87, Issue 1,
2007,  37-56.
%
\bibitem{GGM} Gazzola F.,  Grunau H. C.,  Mitidieri E., Hardy inequalities with
optimal constants and remainder terms, Trans. Amer. Math. Soc. 356,
2004, 2149-2168.
%
%
\bibitem{GK} Ghoussoub N., Kang X.S., Hardy-Sobolev critical elliptic equations with boundary singularities.
Ann. Inst. H. Poincar\'e Anal. Non Lin\'eaire  21  (2004),  no. 6,
767--793.
%
\bibitem{OK}  Opic B. and  Kufner A., "Hardy-type Inequalities", Pitman Research
Notes in Math., Vol. 219, Longman 1990.
%
\bibitem{GT}  Gilbarg D. and  Trudinger N.S., Elliptic partial differential equations of second order.
 $2^{nd}$ edition, Grundlehren 224, Springer, Berlin-Heidelberg-New York-Tokyo (1983).
%
\bibitem{Lind}  Lindqvist P., On the equation $\div(|\n u|^{p-2}\n u)+\l |u|^{p-2}u = 0$, Proc.
Amer. Math. Soc., 109(1) (1990), 157�164. Addendum, ibiden, 116
(2) (1992), 583-584.
%
\bibitem{MMP}  Marcus M.,  Mizel V.J., and  Pinchover Y., Transactions of the American
Mathematical Socity. Volume 350, Number 8, August 1998,  3237-3255.
%
%
\bibitem{NaC}  Nazarov A. I., Hardy-Sobolev Inequalities in a cone, J. Math.
Sciences, 132, (2006), (4), 419-427.
%
\bibitem{Na} Nazarov A.I., Dirichlet and Neumann problems to critical
Emden-Fowler type equations. J Glob Optim (2008) 40, 289-303.
%
\bibitem{PT}{Pinchover Y.,  Tintarev K.},
 { Existence of minimizers for Schr\"{o}dinger operators under domain perturbations with application
  to Hardy's inequality}.
 Indiana Univ. Math. J. 54 (2005), 1061-1074.
%

%
\bibitem{Tib} Tidblom J., A geometrical version of Hardy's inequality for
\r{W}${}^{1,p}(\O) $,
 Proc. Amer. Math. Soc. 132 (2004) 2265-2271.
%
%
%
%
%
%
%
%
%
%
%
\end{thebibliography}
\end{document}